\documentclass[12pt]{amsart}
\usepackage{latexsym,fancyhdr,amssymb,color,amsmath,amsthm,graphicx,listings,comment}
\usepackage[section]{placeins}
 \newtheorem*{satz}{Theorem} \newtheorem{coro}{Corollary}
\definecolor{red1}{rgb}{1,0.9,0.9} \definecolor{blue1}{rgb}{0.9,0.9,1} \definecolor{green1}{rgb}{0.9,1,0.9} 
\definecolor{yellow1}{rgb}{1,1,0.8} \definecolor{yellow2}{rgb}{1,1,0.8}
\def\chapter#1{ \vspace{2mm} \begin{center} \fcolorbox{green1}{green1}{ \parbox{16.2cm}{{\Large {\bf #1}}}} \vspace{2mm} \end{center} }

\def\question#1{ \vspace{2mm} \begin{center} \fcolorbox{blue1}{blue1}{ \parbox{13.2cm}{{\bf Question:} #1}} \vspace{2mm} \end{center} }

\def\satz#1{ \vspace{2mm} \begin{center} \fcolorbox{yellow1}{yellow1}{ \parbox{13.2cm}{{\bf Theorem:} #1}} \vspace{2mm} \end{center} }

\let\paragraph\subsection
 \newcommand{\ZZ}{\mathbb{Z}} \newcommand{\RR}{\mathbb{R}} 
\newcommand{\G}{\mathcal{G}} \newcommand{\R}{\mathcal{R}}  \newcommand{\B}{\mathcal{B}}
\newcommand{\C}{\mathcal{C}} \newcommand{\Sphere}{\mathcal{S}}

\def\Bin#1#2{{#1\choose #2}}

\title{The amazing world of simplicial complexes}
\author{Oliver Knill}
\date{4/22/18, notes for the AMS special session "Discretization in Geometry"}
\address{Department of Mathematics \\ Harvard University \\ Cambridge, MA, 02138 }

\makeindex

\begin{document}
\maketitle

\begin{abstract}
Defined by a single axiom, finite abstract simplicial complexes 
belong to the simplest constructs of mathematics. We look at a
few theorems. 
\end{abstract} 

\chapter{Theorems}

\section{Simplicial complexes}

\paragraph{}
A {\bf finite abstract simplicial complex} $G$ is a finite set of non-empty sets which is 
closed under the process of taking finite non-empty subsets. The {\bf Barycentric refinement} $G_1$
of $G$ is the set of finite subsets of the power set of $G$ which are pairwise 
contained into each other. The new complex $G_1$ defines a finite simple graph $\Gamma=(V,E)$, 
where $V=G$ and $E$ are the pairs where one is contained in the other. 
$G_1$ agrees with the Whitney complex of $\Gamma$, the collection of vertex sets of 
complete sub graphs of $\Gamma$.
\index{Simplicial complex}
\index{Complex ! Simplicial}
\index{Whitney complex}
\index{Complex ! Whitney}
\index{Barycentric refinement}
\index{Refinement ! Barycentric}

\satz{
Barycentric refinements are Whitney complexes.
}

\paragraph{}
Examples of complexes not coming directly from graphs are buildings or matroids.
Oriented matroids are examples of elements of the ring $\R$ generated by simplicial 
complexes. Still, the Barycentric refinement $G_1$ of $G$ always allows to study $G$ with the help
of graph theory.  
\index{Matroid}

\paragraph{}
A subset $H$ of $G$ is called a {\bf sub-complex}, if it is itself a simplicial complex. 
Any subset $H$ {\bf generates} a sub-complex, the smallest simplicial complex in $G$ containing $H$.
The set $\G$ of sub-complexes is a Boolean lattice because it is closed under intersection and
union. The {\bf f-vector} of $G$ is $f=(v_0,v_1, \dots, v_r)$, where $v_k$ is the number of
elements in $G$ with cardinality $k+1$. The integer $r$ is the {\bf maximal dimension} of $G$.
\index{Dimension}
\index{Maximal Dimension}
\index{Generate}
\index{Sub-complex}

\section{Poincar\'e-Hopf}

\paragraph{}
A real-valued function $f:G \to \R$ is {\bf locally injective} if $f(x) \neq f(y)$
for any $x \subset y$ or $y \subset x$. In other words, it is a {\bf coloring} in the graph $\Gamma$ representing
$G_1$. The {\bf unit sphere} $S(x)$ of $x \in G$ is the set $\{ y \in G | (x,y) \in E(\Gamma) \}$.
It is the unit sphere in the metric space $G$, where the distance is the geodesic distance in the
graph representing $G_1$. Define the stable unit sphere $S^-_f(x)= \{ y  \in S(x) \; | \; f(y)<f(x) \}$
and the {\bf index} $i_f(x)=\chi(S^-_f(x))$. The {\bf Poincar\'e-Hopf theorem} is
\index{Index ! Poincar\'e-Hopf}
\index{Locally injective}
\index{Poincar\'e-Hopf}
\index{Theorem ! Poincar\'e-Hopf}
\index{Coloring}

\satz{
$\sum_x i_f(x) = \chi(G)$.
}

\paragraph{}
Classically, for a smooth function with isolated critical points 
on a Riemannian manifold $M$, the same definitions and result apply for
$i_f(x)=\lim_{r \to 0} \chi(S^-_{r,f}(x))$, where $S_r$ is the geodesic
sphere of radius $r$ in $M$ centered at $x$. 
\index{Riemannian manifold}
\index{Manifold ! Riemannian}
\index{Geodesic sphere}

\paragraph{}
If $f(x)={\rm dim}(x)$, then $i_f(x)=\omega(x)$. Poincar\'e-Hopf tells then
that $\chi(G)=\chi(G_1)$. If $f(x)=-{\rm dim}(x)$, then $i_f(x)=\omega(x) (1-\chi(S(x)))$. 
For complexes for which every unit sphere is a $2d$-sphere, we have 
$i_{{\rm dim}}=-i_{-{\rm dim}}$ implying $\chi(G)=0$. 
 
\section{Gauss-Bonnet}

\paragraph{}
Any probability space $\Omega$ of locally injective functions
defines a {\bf curvature} $\kappa(x) = {\rm E}[ i_f(x) ]$. As we have
integrated over $f$, the curvature value $\kappa(x)$ only depends on $x$. 
\index{Probability space}
\index{Integral geometry}
\index{Curvature}

\satz{
$\sum_x \kappa(x) = \chi(G)$.
}

\paragraph{}
If $\Omega$ is the product space $\prod_{x \in G} [-1,1]$ with product
measure so that $f \to f(x)$ are independent identically distributed random variables, 
then $\kappa(x)=K(x)$
is the {\bf Levitt curvature} $1+\sum_{k=0} (-1)^k v_k(S(x))/(k+1)$. The same
applies if the probability space consists of all colorings. 
If $f=1+v_0 t + v_1 t^2 + \dots$ is the generating function of the $f$-vector
of the unit sphere, with anti-derivative $F=t+v_0 t^2/2 + v_1 t^3/3 ...$, 
then $\kappa=F(0)-F(-1)$. Compare $\chi(G)=f(0)-f(-1)$ and 
$\sum_x \chi(S(x)) = f'(0)-f'(-1)$. 
\index{Levitt curvature}
\index{Curvature ! Levitt}

\paragraph{}
If $P$ is the Dirac measure on $f(x)={\rm dim}(x)$, then the curvature is
$\omega(x)$. If $P$ is the Dirac measure on $f(x)=-{\rm dim}(x)$, then 
the curvature is $\omega(x)(1-\chi(S(x))$. 
\index{Dirac measure}

\section{Valuations}

\paragraph{}
A real-valued function $X$ on $\G$ is called a {\bf valuation} if 
$X(A \cap B) + \chi(A \cup B) = \chi(A) + \chi(B)$ for all $A,B \in \G$.
It is called an {\bf invariant valuation} if $X(A)=X(B)$ if $A$ and $B$ are
isomorphic. Let $\G_r$ denote the set of complexes of dimension $r$. 
The {\bf discrete Hadwiger theorem} assures:

\index{Valuation}
\index{Invariant valuation}
\index{Hadwiger theorem}
\index{Theorem ! Hadwiger}

\satz{
Invariant valuations on $\G_r$ have dimension $r+1$. 
}

\paragraph{}
A basis of the space of invariant valuations is given by $v_k: \G \to \RR$. 
Every vector $X=(x_0,\dots, x_r)$ defines a valuation $X(G) = X \cdot f(G)$
on $\G_r$. 
\index{Valuation ! Basis}

\section{The Stirling formula}

\paragraph{}
The $f$-vectors transform linearly under Barycentric refinements. Let ${\rm Stirling}(x,y)$ 
denote the {\bf Stirling numbers} of the second kind. It is the number of times one can
partition a set of $x$ elements into $y$ non-empty subsets. 
The map $f \to S f$ is the {\bf Barycentric refinement operator}

\satz{
$f(G_1)=S f$, where $S(x,y) = {\rm Stirling}(y,x) x!$. 
}

\index{Stirling numbers}
\index{Operator ! Barycentric refinement}
\index{Refinement operator}

\paragraph{}
The matrix is upper triangular with diagonal entries $k!$ the factorial. 
If $X(G) = \langle X,f(G) \rangle = X(G_1)=\langle X,S f(G) \rangle =\langle S^T X f(G)$,
then $X=S^T X$ so that $X$ is an eigenvector to the eigenvalue $1$ of $S^T$. 
The valuation with $X=(1,-1,1,-1, \dots )$ is the {\bf Euler characteristic} $\chi(G)$.
This shows that Euler characteristic is unique: 
\index{Euler characteristic}

\satz{
If $X(1)=1$ and $X(G)=X(G_1)$ for all $G$, then $X=\chi$.
}

\section{The unimodularity theorem}

\paragraph{}
A finite abstract simplicial complex $G$ of $n$ sets defines 
the $n \times n$ {\bf connection matrix} $L(x,y)=1$ if $x \cap y \neq \emptyset$ and $L(x,y)=0$ if 
$x \cap y = \emptyset$. The {\bf unimodularity theorem} is:
\index{Connection matrix}
\index{Theorem ! unimodularity}
\index{Unimodularity theorem}

\satz{
For all $G \in \G$, the matrix $L$ is unimodular. 
}

\section{Wu characteristic}

\paragraph{}
Using the notation $x \sim y$ if $x \cap y \neq \emptyset$, define the {\bf Wu characteristic}
$$ \omega(G) = \sum_{x \sim y} \omega(x) \omega(y) \; . $$
For a complete complex $K_{d+1}$ we have $\omega(K^{d+1})=(-1)^d$. As
every $x \in G$ defines a simplicial complex generated by $\{x\}$, 
the notation $\omega(x)$ is justified. 
\index{Wu characteristic}

\paragraph{}
A complex $G$ is a $d$-complex if 
every unit sphere is a $(d-1)$-sphere. A complex $G$ is a {\bf $d$-complex with boundary}
if every unit sphere $S(x)$ is either $(d-1)$ sphere or a $d-1$-ball. The sets for which $S(x)$
is a ball form the {\bf boundary} of $G$. A complex without boundary is {\bf closed}.
\index{Closed $d$-complex.}
\index{d-complex}
\index{Boundary of a complex}

\satz{
For a $d$-complex $G$ with boundary, $\omega(G) = \chi(G)-\chi(\delta G)$.
}

\paragraph{}
For any $d$ one can define higher {\bf Wu characteristic}
$$ \omega_k(G) = \sum_{x_1 \sim \dots x_k} \omega(x_1) \cdots \omega(x_d) \;  $$
summing over all simultaneously intersecting sets in $G$.

\section{The energy theorem}

\paragraph{}
As $L$ has determinant $1$ or $-1$, the inverse $g=L^{-1}$ is a matrix with integer entries.
The entries $g(x,y)$ are the {\bf potential energy values} between the simplices $x,y$. 
\index{Potential energy}

\satz{
For any complex $G$, we have $\sum_{x} \sum_y g(x,y) = \chi(G)$. 
}

\paragraph{}
This {\bf energy theorem} assures that the total potential energy of a complex is the 
Euler characteristic. 
\index{Energy theorem}
\index{Theorem ! Energy}

\section{Homotopy}

\paragraph{}
The graph $1=K_1$ is {\bf contractible}. Inductively, a graph is 
{\bf contractible} if there exists a vertex $x$ such that both $S(x)$ and $G-x$ are contractible.
The step $G \to G-x$ is a homotopy step. Two graphs are {\bf homotopic} if there exists
a sequence of homotopy steps or inverse steps which brings one into the other. Contractible is
not the same than homotopic to $1$.  A graph $G$ is a {\bf unit ball} if there exists a vertex such that $B(x)=G$.
\index{Contractible}
\index{Homotopic}
\index{Unit ball}

\satz{
If $G$ is a unit ball then it is contractible. 
}

\paragraph{}
It is proved by induction. It is not totally obvious. A {\bf cone extension} $G=D+x$ for the 
{\bf dunce hat} $D$ obtained by attaching a vertex $x$ to $D$ is a ball but we can not take
$x$ away. Any other point $y$ can however be taken away by induction as $G-y$ is a 
ball with less elements.
\index{Dunce hat}
\index{Cone extension}
\index{Collapsible}

\satz{
Contractible graphs have Euler characteristic $1$. 
}

\paragraph{}
The proof is done by induction starting with $G=1$. 
It is not true that the Wu characteristic $\sum_{x \sim y} \omega(x) \omega(y)$ 
is a homotopy invariant as $\omega(K_{n+1})=(-1)^n$.

\section{Spheres}

\paragraph{}
The empty graph $0$ is the $(-1)$ sphere. 
A $d$-sphere $G$ is a $d$-graph for which all $S(x)$ are $(d-1)$ spheres and for which 
there exists a vertex $x$ such that both $G-x$ is contractible. The {\bf 1-skeleton graphs} of the
octahedron and the icosahedron are examples of $2$-spheres. Circular graphs with more than $3$
vertices are $1$-spheres. A simplicial complex $G$ is a {\bf $d$-sphere}, if the graph $G_1$ is a
$d$-sphere. Here is the {\bf polished Euler Gem} 

\index{Empty graph}
\index{Skeleton graph}
\index{Sphere}
\index{Euler gem}
\index{Polished Euler gem}

\satz{
$\chi(G) = 1+(-1)^d$ for a $d$-sphere $G$.
}

\satz{
The join of a $p$-sphere with a $q$-sphere is a $p+q+1$-sphere. 
}

\paragraph{}
The {\bf generating function} of $G$ is $f_G(t)=1+\sum_{k=0}^{\infty} v_k(G) t^{k+1}$ with
$v_k(G)$ being the number of $k$-dimensional sets in $G$. It satisfies
\index{Generating function}

\satz{
$f_{G+H}(t)=f_G(t)+f_H(t)-1$ and $f_{G \oplus H}(t) = f_G(t) f_H(t)$. 
}

For example, for $P_2 \oplus P_2 = S_4$ we have $(1+2t) (1+2t) = 1+4 t + 4t^2$. 

\paragraph{}
Given a $d$-graph. The function ${\rm dim}(x)$ has every point a critical point
and $S^-(x)=\{ y \in S(x) \; | \; f(y)<f(x) \}$ and 
$S^+(x) = \{ y \in S(x) \; | \; f(y)>f(x) \}$ then $S(x) = S^-(x) + S^+(x)$. 

\paragraph{}
Since by definition, a sphere becomes contractible after removing one of its points:

\satz{
$d$-spheres admit functions with exactly two critical points. 
}
\index{Reeb theorem}

Spheres are the $d$-graphs for which the minimal number of critical 
points is $2$. There are no $d$-graphs for which the minimal number of critical
points is $1$. 

\section{Platonic complexes}

\paragraph{}
A {\bf combinatorial CW complex} is an empty or finite ordered sequence of spheres
$G=\{ c_1, \dots, c_n \}$ such that $G_n=\{c_1, \dots, c_n \}$ is obtained 
from $G_{n-1}=\{c_1, \dots, c_{n-1}\}$ by selecting a sphere $c_n$ in $G_{n-1}$
such that $c_n$ is either empty or different from any $c_j$. 
We identify $c_j$ with the cell filling out the sphere. 
Its dimension is 1 plus the dimension of the sphere. 
The Barycentric refinement $G_1$ of $G$ is the Whitney complex 
of the graph with vertex set $G$ and where two vertices $a,b$ are connected 
if one is a sub sphere of the other.
\index{CW complex}
\index{Complex ! CW}

\paragraph{}  
$G$ is a d-sphere if $G_1$ is a $d$-sphere as a simplicial complex. 
A subset $H$ of $G$ is a {\bf sub-complex} of $G$ if $H_1 \subset G_1$ for the refinements. 
\index{Sub complex}

\paragraph{}
The {\bf Levitt curvature} of a cell $c_j$ is $F(0)-F(-1)$, where $F$ is the
anti-derivative of the $f$-generating function
$f=1+t v_0 + t^2 v_1 + \dots$ of the sphere $S(c_j)$. The curvature of a cell 
$x$ in a $2$-sphere is $1-v_0(S(x))/2+v_1(S(x))/3=1-v_0(S(x))/6$. 
The curvature of a cell in a $3$ sphere is $0$. Gauss-Bonnet assures that
the sum of the curvatures is the Euler characteristic. 
\index{Levitt curvature}
\index{Curvature ! Levitt}
\index{Gauss-Bonnet}

\paragraph{}
A $d$-sphere $G$ is called a {\bf Platonic $d$-polytope} if for every $0 \leq k \leq d$
and any cell $c_j$ of dimension $k$, there exists a Platonic $(d-1)$-sphere $P_k$ 
such that $S(x)$ is isomorphic to $P_k$. The $-1$-dimensional sphere $0$ is 
assumed to be Platonic. The $0$-dimensional sphere consisting of two isolated points is Platonic too. The
$1$-dimensional complexes $C_k$ with $k \geq 3$ are the Platonic $1$-spheres. With
$C_3$ one denotes the {\bf $1$-skeleton complex} of $K_3$. Let $P=(p(-1),p(0),p(1),p(2), \dots $ 
denote the number of Platonic $d$-polytopes. In the CW case, we have the familiar classification: 
\index{Platonic ! sphere}
\index{Platonic ! polytope}
\index{Sphere ! Platonic}
\index{Polytope ! Platonic} 
\index{Skeleton complex}

\satz{ ${\rm Platonic}_{CW}=(1,1,\infty,5,6,3,3,3,\dots)$.  }

\paragraph{}
The classification of Platonic polytopes of dimension $d$ which are simplicial complexes is 
easier. There is a unique platonic solid in each dimension except in dimensions
$1,2,3$. In the $1$-dimensional case there are infinitely many. In the
two-dimensional case, only the Octahedron and Icosahedron are Platonic. In the
three dimensional case, there is only the 600 cell and the 16 cells. After that
the curvature condition brings it down to the cross polytopes.

\satz{ ${\rm Platonic}_{SC}=(1,1,\infty,2,2,1,1,1,1,\dots)$. }

\section{Dehn-Sommerville relations}

\paragraph{}
Given a $d$-dimensional complex $G$, any
integer vector $(X_0, \dots, X_d)$ in $Z^{d+1}$ defines
a {\bf valuation} $X(G) = X_0 v_0 + \dots X_d v_d$. By distributing the values $X_k$
attached to each $k$-simplex in $G$ equally to its $k+1$ vertices, we get the {\bf curvature}
$K(x) = \sum_{k=0}^{d} X_k v_{k-1}(S(x))/(k+1)$ for the valuation $X$ and graph $G$ at the
vertex $x$. The formula $\sum_{x \in V} K(x) =X(G)$ is the {\bf Gauss-Bonnet theorem} for $X$.

\paragraph{}
In the case $X(G)=v_1$, the curvature is the vertex degree divided by $2$ and the formula reduces
to the ``Euler handshake". If $X=v_d$ is the volume of $G$, then $K(x)$ is the number
of $d$-simplices attached to $x$ divided by $d+1$. In the case $X=(1,-1,1,-1,\dots)$,
$X$ is the {\bf Euler characteristic} and $K$ is the discrete analogue of the Euler form
in differential geometry entering the Gauss-Bonnet-Chern theorem. For $d$-graphs, there are
some valuations which are zero. Define the {\bf Dehn-Sommerville} valuations
$$ X_{k,d} = \sum_{j=k}^{d-1} (-1)^{j+d} \Bin{j+1}{k+1} v_j(G) + v_k(G)  \; . $$

\index{Euler handshake}
\index{Gauss-Bonnet ! valuation} 
\index{Valuation}
\index{Dehn-Sommerville valuation} 

\satz{
For $d$-graphs, the Dehn-Sommerville curvatures are zero. 
}

\paragraph{}
The proof is by noticing that the curvature of $X_{k,d}$ is
$K(x) = X_{k-1,d-1}(S(x))$. This follows from the relation
$$   X_{k+1,d+1}(l+1)/(l+1) = X(k,d)(l)/(k+2)  \; .  $$
Use Gauss-Bonnet and induction using the fact that the unit sphere of a geometric
graph is geometric and that for $d=1$, a geometric graph is a cyclic graph $C_n$
with $n \geq 4$. For such a graph, the Dehn-Sommerville valuations are zero.

\section{Dual Connection matrix}

\paragraph{}
Define the {\bf dual connection matrix} $\overline{L}(x,y)=1-L(x,y)$. 
It is the adjacency graph of a {\bf dual connection graph}, where two 
simplices are connected, if they do not intersect. 
\index{Dual | connection matrix}
\index{Connection matrix ! dual}
\index{Connection graph ! dual}

\satz{
$1-\chi(G) ={\rm det}(-L \overline{L})$.
}

\paragraph{}
Let $E$ be the constant matrix $E(x,y)=1$. The result follows from 
unimodularity ${\rm det}(L)={\rm det}(g)$ and the energy theorem telling that
$\overline{L} g =(E-L) g = E g - 1$ has the eigenvalues
of $Eg$ minus $1$ which are $\chi(G)$ and $0$. Assume $G$ has $n$ sets: 

\satz{
$-L \overline{L}$ has $n-1$ eigenvalues $1$ and one eigenvalue $1-\chi(G)$.
}

\paragraph{}
The above formula is not the first one giving the Euler characteristic as a determinant
of a Laplacian. \cite{BeraKMukherjee} show, using a formula of Stanley, that if $A(x,y)=1$ 
if $x$ is not a subset of $y$ and $A(x,y)=0$ else, then $1-\chi(G) = \det(A)$.

\section{Alexander Duality}

\paragraph{}
The {\bf Alexander dual} of $G$ is the simplicial complex
$G^* = \{ x \subset V \; | \; x^c \notin G \}$. It is the complex
generated by the complements $x^c$ of the sets $x$ in $G$. For the complete
complex $K_d$, the dual is the empty complex. In full generality one has
for the Betti numbers $b_k(G)$
\index{Alexander dual}
\index{Dual  ! Alexander}
\index{Betti number ! dual}

\satz{$b_k(G)=b_{n-3-k}(G)$, $k=1, \dots, n-1$}

\paragraph{}
In order to have content, this needs $n \geq 5$. It works for $G=C_5$ already, 
where $G^*$ is the complement of a circle in a $3$-sphere. 
The combinatorial Alexander duality is due to Kalai and Stanley. 
\index{Kalai}
\index{Stanley}

\section{Sard}

\paragraph{}
Given a locally injective function $f$ on a graph $G=(V,E)$, define for $c \notin f(V)$ 
the {\bf level surface} $\{ f=c \}$ as the subgraph of the Barycentric refinement of 
$G$ generated by simplices $x$ on which $f$ changes sign. Remember that $G$ is a $d$-graph
if every unit sphere $S(x)$ is a $(d-1)$-sphere. A discrete Sard theorem is:
\index{Theorem ! Sard}
\index{Level surface}
\index{Contour surface}
\index{Sard Theorem}

\satz{
For a d-graph, every level surface is a $(d-1)$-graph.
}

If $G$ is a finite abstract simplicial complex, then $f: G \to R$ defines a 
function on the Barycentric refinement $G_1$ and the level surface is defined like that. 
This result has practical value as we can define discrete versions of classical surfaces. 

\paragraph{}
Given a finite set of functions $f_1, \dots, f_k$ on the vertex sets of Barycentric refinements 
$G_1, \dots, G_k$ of a simplicial complex, we can now look at the $(d-k)$-graph
$\{ f=c \} = \{ f_1=c_1, \dots, f_k=c_k \}$. Unlike in the continuum case, where the result only holds
for almost all $c$, this holds for all $c$ disjoint from the range. 

\satz{
Given $f_j:G_j \to \R, j \leq k$, then $\{ f=c \}$ is a $(d-k)$-complex.
}

\section{Bonnet and Synge}

\paragraph{}
The topic of positive curvature complexes is analog to the continuum. Still, 
it would be nice to have entirely combinatorial proofs of the results in the
continuum. 

\paragraph{}
Let $G$ be a $d$-complex so that every unit sphere is a $(d-1)$ sphere. 
A {\bf geodesic $2$-surface} is a subcomplex if the embedded graph does not contain a 3-simplex.
$G$ has positive sectional curvature if every geodesic embedded wheel graph $W(x)$ has
interior curvature $\geq 5/6$.  The {\bf geomag lemma} is
that any wheel graph in a positive curvature $G$ can be extended to an embedded 2-sphere.
\index{Geomag lemma}
\index{Positive curvature}
\index{Sectional curvature}

\paragraph{}
An elementary analog of the {\bf Bonnet  theorem}
\index{Bonnet theorem}

\satz{A positive curvature complex has diameter $\leq 4$.}

\paragraph{}
The simplest analog of {\bf Synge theorem} is 

\satz{A positive curvature complex is simply connected.}

\paragraph{}
The reason for both statements is the {\bf geomag lemma} stating 
that any closed geodesic curve can be extended to a 2-complex which is a sphere and so simply connected.
The strict curvature assumption as we can not realize a projective plane yet with so few cells. With
weaker assumptions getting closer to the continuum, we also have to work harder: 

\paragraph{}
Define more generally the {\bf sectional curvature to be $\geq \kappa$} if there exists $M$ such that
the total interior curvature of any geodesic embedded $2$-disk with $M$ interior points is 
$\geq \delta \cdot M$ and such that every geodesic embedded wheel graph $W(x)$ has non-negative
interior curvature. A complex has {\bf positive curvature} if there exists $\kappa>0$ such
that $G$ has sectional curvature $\geq \kappa$. The maximal $\kappa$ which is possible is 
then the "sectional curvature bound". 

\paragraph{}
An embedded $2$-surface of positive sectional curvature $\kappa$ must then have surface area 
$\leq 2/\kappa$. The classical theorem of Bonnet assures that a Riemannian manifold of 
positive sectional curvature is compact and satisfies an upper diameter bound $\pi/\sqrt{k}$. 
An analog bound $C/\sqrt{k}$ should work in the discrete. 

\paragraph{}
Having a notion of sectional curvature allows to define {\bf Ricci curvature} of an edge $e$
as the average over all sectional curvatures over all wheel graphs passing through $e$. The 
{\bf scalar curvature} at a vertex $x$ is the average Ricci curvatures over all edges $e$ 
containing $x$. The {\bf Hilbert functional} is then the total scalar curvature. Unlike in 
Regge calculus, all these notions are combinatorial and do not depend on an embedding. 
\index{Ricci curvature}
\index{Scalar curvature}
\index{HIlbert functional}

\section{An inverse spectral result}

\paragraph{}
Let $p(G)$ denote the number of positive eigenvalues of the connection 
Laplacian $L$ and let $n(G)$ the number of negative eigenvalues of $L$.
\index{Eigenvalues ! Positive}
\index{Eigenvalues ! Negative}

\satz{
For all $G \in \G$ we have $\chi(G) = p(G)-n(G)$.
}

\index{Spectral formula}
\index{Hearing ! Euler characteristic} 

\paragraph{}
The proof checks this by deforming $L$ when adding a new cell.
This result implies that Euler characteristic is a logarithmic potential energy 
of the origin with respect to the spectrum of $i L$. 
\index{Logarithmic energy}
\index{Logarithmic potential}

\satz{
$\chi(G)={\rm tr}(\log(i L)) (2 \pi/i)$. 
}

\paragraph{}
The proof shows also that after a $CW$ ordering of the sets in a finite abstract
simplicial complex, one can assign to every simplex a specific eigenvalue and so
eigenvector of $L$. 

\section{The Green star formula}

\paragraph{}
Given a simplex $x \in G$, the {\bf stable manifold} of the
dimension functional ${\rm dim}(x)$ is $W^-(x)=\{ y \in G \; | \ ; y \subset x \}$. 
The {\bf unstable manifold} $W^+(x) = \{ y \in G \; | \; x \subset y \}$ is known as the
{\bf star} of $x$. Unlike $W^-(x)$ which is always a simplicial complex, the star
$W^+(x)$ is in general not a sub complex of $G$. 
\index{Stable manifold}
\index{Unstable manifold}
\index{Star}

\satz{
$g(x,y) = \omega(x) \omega(x) \chi(W^+(x) \cap W^+(y))$. 
}

\index{Green star formula}

\paragraph{}
In comparison, we have $W^-(x) \cap W^+(x)=\omega(x)$ and  
$L(x,y)=\chi(W^-(x) \cap W^-(y))$. The to $L$ similar matrix 
$M(x,y)=\omega(x) \omega(x) \chi(W^-(x) \cap W^-(y))$
satisfies $\sum_x \sum_y M(x,y)=\omega(G)$, the Wu characteristic. 

\section{Wu characteristic}

\paragraph{}
The {\bf Euler characteristic} $\chi(G)=\omega_1(G) = \sum_{x \in G} \omega(x)$ of $G$ is the simplest of a sequence of
combinatorial invariants $\omega_k(G)$. The second one, $\omega(G) = \sum_{x,y, L(x,y)=1} \omega(x) \omega(y)$, is
the {\bf Wu characteristic} of $G$. The valuation $\chi$ is an example of a linear valuation, while $\omega$ is a
{\bf quadratic valuation}.
The Wu characteristic also defines an {\bf intersection number} $\omega(A,B)$ between sub-complexes.
\index{Wu characteristic}
\index{Intersection number}
\index{Valuation ! quadratic}

\paragraph{}
All multi-linear valuations feature Gauss-Bonnet and Poincar\'e-Hopf theorems,
where the curvature of Gauss-Bonnet is an index averaging. For example, with
$K(v) = \sum_{v \in x, x \sim y} \omega(x) \omega(y)/(|x|+1)$
The Gauss-Bonnet theorem for Wu characteristic is
\index{Multi-linear valuation}
\index{Wu ! curvature}
\index{Curvature ! Wu}
\index{Gauss-Bonnet ! multi-linear valuation}

\satz{
$\omega(G) = \sum_v K(v)$.
}

\index{Gauss-Bonnet ! Wu characteristic}

\section{The boundary formula}

\paragraph{}
We think of the {\bf internal energy} $E(G)=\chi(G)-\omega(G)$ as a sum of
{\bf potential energy} and {\bf kinetic energy}. 
A {\bf d-complex} is a simplicial complex $G$
for which every $S(x)$ is a $(d-1)$-sphere. A {\bf d-complex with boundary} 
is a complex $S(x)$ is either a $(d-1)$-sphere or a $(d-1)$-ball for every $x \in G$. 
\index{Complex with boundary}
\index{Boundary}
\index{Complex ! d-complex}

\paragraph{}
The $d$-complexes are {\bf discrete $d$-manifolds} and $d$-complexes with boundary is a
discrete version of a {\bf $d$-manifold with boundary}. We denote by $\delta G$ the 
{\bf boundary} of $G$. It is the $d-1$ complex consisting of boundary points. 
By definition, $\delta \delta G=0$, the empty complex. The reason is that the 
boundary of a complex is closed, has no boundary. We can reformulate the formula
given below as 

\index{Internal energy}
\index{Energy ! internal}
\index{Energy ! potential}
\index{Energy ! kinetic}
\index{Discrete manifold with boundary}

\satz{
If $G$ is a $d$-complex with boundary then $E(G)=\chi(\delta(G))$.
}

\index{Complex}

\index{Boundary formula}

\paragraph{}
If $G$ is a $d$-ball, then $\delta G$ is a $(d-1)$-sphere
and $E(G)=1+(-1)^{d-1}$, by the {\bf polished Euler gem formula}. 
\index{Ball}
\index{Polished Euler gem}
\index{Euler gem formula}

\section{Zeta function}

\paragraph{}
For a one-dimensional complex $G$, there is a {\bf spectral symmetry} 
which will lead to a {\bf functional equation}:
\index{Spectral symmetry}
\index{Functional equation}

\satz{
If ${\rm dim}(G)=1$, then $\sigma(L^2)=\sigma(L^{-2})$. 
}

\paragraph{}
If $H$ is a Laplacian operator with non-negative spectrum 
like the {\bf Hodge operator} $H$ or {\bf connection operator} $L$,
one can look at its {\bf zeta function}
$$ \zeta_H(s) = \sum_{\lambda \neq 0} \lambda^{-s} \; , $$
where the sum is over all non-zero eigenvalues of $H$ or $L^2$. In
the connection case, we take $L^2$ to have all eigenvalues positive.
\index{Zeta function}
\index{Spectral zeta function}

\paragraph{}
The case of the connection Laplacian is especially interesting because one does not have
to exclude any zero eigenvalue. The {\bf connection zeta function} of $G$ is defined as
$\zeta(s) = \sum_{\lambda} \lambda^{-s}$, where the sum is 
over all eigenvalues $\lambda$ of $L^2$. It is an entire function in $s$. 
\index{Zeta function ! connection}
\index{Entire function}
\index{Connection zeta function}

\satz{
If ${\rm dim}(G)=1$, then $\zeta(s)=\zeta(-s)$.
}

\index{Functional equation ! Spectral zeta}

\paragraph{}
When doing Barycentric refinement steps, the zeta function converges to an explicit function.
$$ \zeta(it) = \int_0^1 \frac{2 \cos \left(2 t \log \left(\sqrt{4 v^2+1}+2 v\right)\right)}{\pi  \sqrt{1-v} \sqrt{v}} \, dv \;  . $$
It is a {\bf hypergeometric series} $\zeta(2s)=\pi \, _4F_3\left(\frac{1}{4},\frac{3}{4},-s,s;\frac{1}{2},\frac{1}{2},1;-4\right)$. 
\index{Hypergeometric series}
\index{Limiting zeta function}

\section{The Hydrogen formula}

\paragraph{}
Given a simplicial complex $G$, let $\Lambda_k(G)$ denote the set of
real valued functions on $k$-dimensional simplices. It is a $v_k$-dimensional
vector space. Define the $v_k \times v_{k+1}$ matrices $d_k(x,y)=1$ if $x \subset y$
and $d_k(x,y)=0$ else. It is the {\bf sign-less incidence matrix}. It can be extended to 
a $n \times n$ matrix $d$ so that $d=d_0+d_1 + \cdots + d_r$ and $D=d+d^*$ and $H=(d+d^*)^2$,
the {\bf sign-less Dirac} and {\bf sign-less Hodge operator}. In the one-dimensional case, 
we have $H=d^* d + d d^*$. The {\bf Hydrogen relations} are
\index{Sign-less incidence matrix}
\index{Incidence matrix ! sign-less}
\index{Hydrogen relation}

\satz{
If ${\rm dim}(G)=1$, then $L-L^{-1}=H$.
}

\paragraph{}
The relation allows to relate the spectra of $L$ and $H$. It allows to estimate the 
spectral radius or give explicit formulas for the spectrum of the connection Laplacian in 
the circular case. This is needed to get the explicit {\bf dyadic zeta function}
\index{Dyadic zeta function}
\index{Zeta function ! Dyadic}

\paragraph{}
Let $S(x)$ denote the {\bf unit sphere} of a simplex $x \in G$. While $S(x)$ is at first a subset of $G$, 
it generates a sub-complex in $G_1$. As $g(x,x) = 1-\chi(S(x)) = \chi(W^+(x))$, we have a functional 
$\sum_x \chi(S(x))$ of Dehn-Sommerville type. With $f(t) = 1+\sum_{k=1}^{\infty} v_{k-1} t^k
=1+v_0 t + v_1 t^2+ v_2 t^3 + \cdots$, the Euler characteristic of $G_1$ can be written as 
$\chi(G)=f(0)-f(-1)$. The following result holds for any simplicial complex: 
\index{Unit sphere}

\satz{
${\rm tr}(L-L^{-1}) = \sum_{x} \chi(S(x)) = f'(0)-f'(-1)$. 
}

\paragraph{}
Compare that the Levitt curvature at a point $x$ was $F(0)-F(-1)$, where $F$ is the anti-derivative
of the generating function of $S(x)$. 
\index{Levitt curvature}
\index{Curvature ! Levitt}

\section{Brouwer-Lefschetz}

\paragraph{}
The {\bf exterior derivative} $d$ for $G$ defines the {\bf Dirac operator} $D=d+d^*$ of $d$. The Hodge
Laplacian $H=D^2$ splits into a direct sum $H_0 \oplus H_1 \cdots H_d$. 
The null space of $H_k$ is isomorphic to the $k$'th cohomology group $H^k(G)={\rm ker}(d_k)/{\rm im}(d_{k-1})$. 
Its dimension $b_k$ is the $k$'th Betti number. The {\bf Euler-Poincar\'e} relation assures that the 
cohomological Euler characteristic $\sum_k (-1)^k b_k$ is equal to the Euler characteristic. 
\index{Incidence matrix}
\index{Exterior derivative}
\index{Euler-Poincar\'e}
\index{Cohomological Euler characteristic}

\paragraph{}
An {\bf endomorphism} $T$ of $G$ is a map from $G$ to $G$ which preserves the order structure. 
It is an automorphism if it is bijective. An endomorphism $T$ induces a linear map on 
cohomology $H^k(G)$. The super trace of this map is 
the {\bf Lefschetz number} $\chi(T,G)$ of $T$. Given a fixed point $x \in G$ of $T$, 
its {\bf Brouwer index} is defined as $i_T(x)=\omega(x) {\rm sign}(T|x)$. Now
\index{Endomorphism}
\index{Automorphism}
\index{Brouwer index}
\index{Lefschetz number}
\index{index ! Brouwer}
\index{Theorem ! Lefschetz}
\index{Lefschetz fixed point theorem}
\index{Fixed point theorem ! Lefschetz}
\index{Fixed point theorem ! Brouwer}
\index{Brouwer fixed point theorem}

\satz{
$\chi(T,G) = \sum_{x=T(x)} i_T(x)$. 
}

\paragraph{}
A special case is $T=1$, where $\chi(1,G)=\chi(G)$ and $i_T(x)=\omega(x)$. 
The Brouwer-Lefschetz fixed point theorem is then the Euler-Poincar\'e theorem. 
\index{Euler-Poincar\'e theorem}

\section{McKean-Singer}

\paragraph{}
The {\bf super trace} ${\rm str}(A)$ of a $n \times n$ matrix defined for a complex with 
$n$ sets is defined as $\sum_{x \in G} \omega(x) L(x,x)$. By definition, 
we have ${\rm str}(1)={\rm str}(L)$. For the Hodge operator $H=D^2=(d+d^*)^2$ we have the
{\bf McKean-Singer formula}: 
\index{Super trace}
\index{McKean-Singer}
\index{Theorem ! McKean-Singer}

\satz{
${\rm str}(\exp(-t H)) = \chi(G)$ for all $t$.
}

\paragraph{}
The reason is that ${\rm str}(H^k)=0$ for $k > 0$, implying ${\rm str}(\exp(t H))={\rm str}(1)=\chi(G)$. 
The McKean-Singer identity is very important as it allows to give almost immediate {\bf proofs} of 
the Lefschetz formulas in any framework in which the identity holds. We proposed in \cite{DiscreteAtiyahSingerBott}
to define a discrete version of a {\bf differential complex} as McKean-Singer enables Atiyah-Singer or 
Atiyah-Bott like extensions of Gauss-Bonnet or Lefschetz. They are caricatures of the heavy theorems in the continuum. 
\index{Differential complex}
\index{Atiyah-Bott}
\index{Atiyah-Singer}

\paragraph{}
The Hodge operator $H=(d+d^*)^2$ and the connection operator $L$ live on the same 
finite dimensional Hilbert space. There is no cohomology associated to $L$. But 
for the connection operator $L$, there is still a localized version of McKean-Singer:

\satz{
${\rm str}(L^k) = \chi(G)$ for $k =-1,0,1$. 
}

\index{McKean-Singer formula ! connection}

\section{Barycentric limit}

\paragraph{}
A matrix $L$ with eigenvalues $\lambda_0 \leq \lambda_2 \leq \cdots \leq \lambda_{n-1}$ 
defines a {\bf spectral function} $F(x)=\lambda_{[nx]}$ on $[0,1)$, where $[t]$ is the 
floor function giving the largest integer smaller or equal than $t$. The inverse function $k(x)=F^{-1}(x)$
is called the {\bf integrated density of states} of $L$ and $\mu=k'$ is the {\bf density of states}.
The sequence $G_k$ of Barycentric refinements of $G$ defines a sequence of operators
$L_k$ and so a sequence of spectral functions $F_n(x)$. Let $\G_r$ denote the set of complexes
of dimension $r$. The following {\bf spectral universality} is a {\bf central limit theorem}: 
\index{Spectral function}
\index{Density of states}
\index{Integrated density of states}
\index{Spectral universality}

\satz{
$\exists F=F(r)$ such that $F_n(G)  \to_{L^1} F$ for all $G \in \G_r$.
}

\paragraph{}
For $r=1$, we know  $F(x)=4 \sin^2(\pi x/2)$. The function is important as 
it conjugates the {\bf Ulam map} $z \to 4x(1-x)$ to a linear function $T(F(x))=F(2x)$.
The measure $\mu$ maximizes {\bf metric entropy} of the Ulam map and is equal to
the {\bf topological entropy} which is $\log(2)$ for $T$.
\index{Metric entropy}
\index{Ulam map}
\index{Topological entropy}

\paragraph{}
We think of $G_n \to G_{n+1}$ as a renormalization step like adding and normalizing two
independent random variables. The result can be seen as a {\bf central limit theorem}.
\index{Central limit theorem}

\section{The join monoid}

\paragraph{}
The {\bf join} $G + H$ of two complexes $G,H$ is the complex $G \cup H \cup \{ x\cup y, x \in H, y \in G\}$. 
For graphs it is known as the {\bf Zykov sum}. 
Given graphs $G=(V,E),H=(W,F)$ then the sum is $(V \cup W,E \cup F \cup \{ (a,b) \; | \; a \in V, b \in W\})$.
If $\overline{G}$ denotes the {\bf complement graph} and $+$ the disjoint union,
then $\overline{G \oplus H} = \overline{G} + \overline{H}$. 
\index{Complement graph}

\paragraph{}
The {\bf join} of two simplicial complexes $G,H$ is defined as
the complex generated by $G + H  = G \cup H \cup \{ x \cup y \; | \; x \in G, y \in H\}$.
Let $f_G(t) = 1+v_0 t + v_1 t^2 + \dots $ denote the 
{\bf generating function} of $G$: then we have the multiplication formula: 
\index{Join monoid}
\index{Monoid ! Zykov}
\index{Monoid ! join}
\index{Zykov sum}

\satz{
$f_{G + H}(t) = f_G(t) f_H(t)$.
}

\paragraph{}
This gives $1-\chi(G)=f_G(-1)$. The dimension function on $G$
not only defines a coloring on $G_1$, it also defines a hyperbolic splitting of the unit spheres. Let 
$S^-(x)=\{ y \in S(x), {\rm dim}(y)<{\rm dim}(x) \}$ and 
$S^+(x)=\{ y \in S(x), {\rm dim}(y)>{\rm dim}(x) \}$. We call them the {\bf stable sphere} and {\bf unstable sphere}. 
\index{stable sphere}
\index{unstable sphere}

\satz{
$S(x) = S^-(x) +  S^+(x)$.
}

\index{Hyperbolic splitting}

\paragraph{}
It follows that 
$g(x,x)=1-\chi(S(x)) = (1-\chi(S^-(x)) ) (1-\chi(S^+(x)) ) = \omega(x) (1-\chi(S^+(x)))$. 
This implies that ${\rm str}(L^{-1}) = \sum_x (1-\chi(S^+(x))) = \chi(G)$ because 
this is the sum over the Poincar\'e-Hopf indices of the function $-{\rm dim}$. 

\paragraph{}
The join monoid is isomorphic to the additive monoid of disjoint union. The {\bf zero element}
is $0$, the $-1$ sphere. 
One can show by induction that if $H$ is contractible and $K$ arbitrary 
then $H+K$ is contractible. This implies:
\index{Sphere monoid}
\index{Monoid ! sphere}

\satz{The join $G$ of two spheres $H+K$ is a sphere.}

\paragraph{}
For example, the join of two zero dimensional spheres $P_2$ is the circle $C_4$. 
The join of two circles a three sphere. It is not the dimension but the {\bf clique number} 
${\rm dim}(G)+1$ which is additive. The clique number of the $-1$ sphere $0$ is $0$.
\index{Clique number}  

\section{The strong ring}

\paragraph{}
The {\bf addition} $A+B$ of two complexes is the disjoint union. 
The empty complex $0$ is the {\bf zero element}. The {\bf Cartesian product}
$G \times H$ is not a simplicial complex any more. We can look at the 
ring $\R$ generated by simplicial complexes. It has the one point complex $1=K_1$
as {\bf one element}. Connected elements are the {\bf additive primes}, simplicial complexes are
{\bf multiplicative primes}. The {\bf Hodge operator} $H$ and the 
{\bf connection operator} $L$ can both be extended to the ring $\R$.

\index{Addition ! complex}
\index{Signed complex}
\index{Cartesian product}
\index{Multiplicative primes}
\index{Primes ! Additive}
\index{Primes ! Multiplicative}
\index{Strong ring}
\index{One element} 
\index{Ring ! strong}

\satz{
$\sigma(H(A \times B)) = \sigma(H(A)) + \sigma(H(B))$,
}

\paragraph{}
Furthermore: 

\satz{
$\sigma(L(A \times B)) = \sigma(L(A)) \cdot \sigma(L(B))$
}

\section{Kuenneth formula}

\paragraph{}
The {\bf Betti numbers} of a signed complex $b_k(G)$ are now signed with $b_k(-G)=-b_k(G)$.
The maps assigning to $G$ its Poincar\'e polynomial $p_G(t)=\sum_{k=0} b_k(G) t^k$ or
{\bf Euler polynomial} $e_G(t)=\sum_{k=0} v_k(G) t^k$ are ring homomorphisms from $R$ to $\ZZ[t]$.
Also $G \to \chi(G)=p(-1)=e(-1) \in \ZZ$ is a ring homomorphism.
\index{Betti numbers ! Signed}
\index{Signed Betti numbers}
\index{Poincar\'e polynomial}
\index{Euler polynomial}
\index{Polynomial ! Euler}
\index{Polynomial ! Poincar\'e}

\satz{ $e_G$ and $p_G$ are ring homomorphisms $\R \to \ZZ[t]$.  }

\paragraph{}
The {\bf Kuenneth formula} for cohomology groups $H^k(G)$ is explicit via Hodge: 
a basis for $H^k(A \times B)$ is obtained from a basis of the factors.
The product in $R$ produces the strong product for the connection graphs. These relations generalize to
Wu characteristic. $R$ is a subring of the full {\bf Stanley-Reisner ring} $S$, a
subring of a quotient ring of the polynomial ring $Z[x_1,x_2, \dots ]$. An object $G \in R$ can
be visualized by ts Barycentric refinement $G_1$ and its connection graph $G'$.
\index{Kuenneth formula}
\index{Formula ! Kuenneth}
\index{Stanley Reisner ring}
\index{Wu characteristic}

\paragraph{}
Theorems like Gauss-Bonnet, Poincar\'e-Hopf or Brouwer-Lefschetz for Euler and Wu characteristic
extend to the strong ring. The isomorphism $G \to G'$ to a subring of the strong Sabidussi ring shows that
the multiplicative primes in $R$ are the simplicial complexes and
that connected elements in $R$ have a unique prime factorization.

\paragraph{}
The {\bf Sabidussi ring} is dual to the Zykov ring. The Zykov join was the addition which is a
sphere-preserving operation. The Barycentric limit theorem implies that the connection 
Laplacian remains invertible in the limit. 
\index{Sabidussi ring}
\index{Ring ! Sabidussi}

\section{Dimension}

\paragraph{}
The {\bf inductive dimension} of a graph is defined inductively as 
${\rm dim}(G) = 1+\sum_{v \in V} {\rm dim}(S(x))/|V|$. For a general complex
$G$ we can define ${\rm dim}(G) = {\rm dim}( G_1)$, as $G_1$ is now
the Whitney complex of a graph. We have ${\rm dim}(G) \leq {\rm maxdim}(G)=
\max_{x \in G} (|x|-1)$, where the right hand side is the {\bf maximal
dimension}.
\index{Dimension ! inductive}
\index{Dimension ! maximal}
\index{Inductive dimension}
\index{Maximal dimension}

\satz{
${\rm dim}(A \times B) = {\rm dim}(A) + {\rm dim}(B)$. 
}

\paragraph{}
Under Barycentric refinements, the inductive dimension can only increase.

\satz{
${\rm dim}(G_1) \geq {\rm dim}(G)$
}

\index{Inductive dimension ! monotonicity}

\paragraph{}
The reason is that higher dimensional complexes have more off-springs than
smaller dimensional ones. 

\paragraph{}
This implies a inequality which resembles the corresponding inequality for 
{\bf Hausdorff dimension} in the continuum: 

\satz{
${\rm dim}( (A \times B)_1 )  \geq {\rm dim}(A) + {\rm dim}(B)$.
}

\index{Hausdorff dimension}

\section{Random complexes}

\paragraph{}
Given a probability space of complexes, one can study the expectations
of random variables. The simplest probability space is the 
{\bf Erd\"os-R\'enyi space} $E(n,p)$ of random graphs equipped with the Whitney complex. 
Define the polynomials $d_n(p)$ of degree $\Bin{n}{2}$ as
$$ d_{n+1}(p) = 1+\sum_{k=0}^n \Bin{n}{k} p^k (1-p)^{n-k} d_k(p) \; , $$
where $d_0=-1$. We can now estimate the inductive dimension.
\index{Erdoes-Renyi space}
\index{Random graphs}

\satz{
${\rm E}_{G(n,p)}[{\rm dim}]) = d_n(p)$. 
}

\paragraph{}
As the Euler characteristic is one of the most important functionals, we want to
estimate its expectation:

\satz{
$$ {\rm E}_{G(n,p)}[\chi] = \sum_{k=1}^n (-1)^{k+1} \Bin{n}{k} p^{\Bin{k}{2}} \; . $$
}

\index{Expectation of Euler characteristic}

\paragraph{}
We don't yet know the expectation value of the Wu characteristic on $E(n,p)$. 

\section{Lusternik-Schnirelmann} 

\paragraph{}
A complex $G$ is {\bf contractible} if there exists $x \in G$ such that both the unit sphere $S(x)$ 
as well as the complex $G \setminus x$ are contractible. A complex is {\bf homotopic to K=1} if there there exists
a complex $H$ such that $H$ is contractible to both $G$ and $K$. The {\bf dunce hat} is 
an example of a complex homotopic to $1$ which is not contractible. The minimal number
of contractible subcomplexes of $G$ covering $G$ is called the Lusternik-Schnirelman category
of $G$. 
\index{Contractible}
\index{Homotopic}

\paragraph{}
A $x \in G$ is called a {\bf critical point} of a function $f$ if $S^-_f(x)$ is not
contractible. The minimal number of critical points which a function $f$ on $G$ can have is 
denoted by ${\rm cri}(G)$. 

\paragraph{}There is a graded multiplication $H^k(G) \times H^l(G)  \to H^{k+l}(G)$ 
called the {\bf cup product}. If $m-1$ is the maximal number of $p>0$-forms $f_1, \dots, f_{m-1}$
for which $f_1 \cup \cdots \cup f_{m-1}$ is not zero, then $m$ is called the {\bf cup length} 
of $G$. 
\index{Cup product}
\index{Theorem ! Lusternik-Schnirelmann}

\paragraph{}
The following result, established with Josellis in 2012 is completely analog to the 
continuum. 

\satz{
${\rm cup}(G) \leq {\rm cat}(G) \leq {\rm cri}(G)$.  
}

\paragraph{}
For any critical point $x_i$, we can form the maximal complex $G_i$ which does not contain an other
critical point. Each $U_i$ is contractible and cover $G$. This proves ${\rm cat}(G) \leq {\rm cri}(G)$.
If ${\rm cat}(G)=n$, let $\{U_k \; \}_{k=1}^n$ be a {\bf Lusternik-Schnirelmann cover}.
Given a collection of $k_j \geq 1$-forms $f_j$ with $f_1 \wedge f_2 \dots \wedge \dots f_n \neq 0$.
Using coboundaries we can achieve that for any simplex $y_k \in U_k$, we can change $f$ in the
same cohomology class $f$ so that
$f(y_k)=0$. Because $U_k$ are contractible in $G$, we can render $f$ zero in $U_k$. This shows that
we can choose $f_k$ in the relative cohomology groups $H^k(G,U_k)$ meaning that we can
find representatives $k_j$ forms $f_j$ which are zero on each $p_{k_j}$ simplices
in the in $G$ contractible sets $U_k$. But now, taking these
representatives, we see $f_1 \wedge \cdots \wedge f_n = 0$. This shows ${\rm cup}(G) \leq n$.
\index{Cover ! Lusternik-Schnirelmann}

\section{Morse inequality}

\paragraph{}
A locally injective scalar function $f$ on the vertex set of a d-graph is called a 
{\bf Morse function}, if $S^-_f(x)$ is a sphere for every $x$. 
The {\bf Morse index} is $m(x)=1+{\rm dim}(S_f^-(x))$. The {\bf Poincar\'e-Hopf index} is $(-1)^{m(x)}$. 
For example, if $d=2$, and $S_f^-(x)$ is $0$-dimensional, then $m(x)=1$ and $i_f(x)=-1$. 
A function $f$ on an abstract simplicial $d$-complex $G$ is a Morse function if it is a Morse
function on the graph $G_1$. 
\index{Morse function}
\index{Morse index}
\index{Poincare-Hopf index}

\satz{Every $d$-complex admits a Morse function.}

\paragraph{}
We can build up $G$ as a {\bf discrete $CW$-complex}. The number at which a simplex $x$ has been added 
is a Morse function as $S(x)$ and $S^-(x)$ are both spheres. 
Also the function ${\rm dim}(x)$ is a Morse function. For $d$-complexes, the stars of two simplices
intersect in a simplex so that: 
\index{CW complex}

\satz{For a $d$-complex, the Green function takes values $1,-1,0$.}

We have $g(x,y) = \omega(x) \omega(y) \chi(W^+(x) \cap W^+(y))$. We have 
$W^+(x) \cap W^+(y) = (1-S^+(x))(1-\chi(S^+(x))$ which is in $\{-1,1\}$ if there
is an intersection and $0$ if not. 
Let $b_k(G)$ denote the $k$'th Betti number. Let $c_k(G)$ denote the number
of critical points of index $k$. Here are the {\bf weak Morse inequalities}:

\satz{$b_k(G) \leq c_k(G)$.}

We even have the {\bf strong Morse inequalities}
\index{Strong Morse inequalities}
\index{Morse inequalities ! Strong}

\satz{ $(-1)^p \sum_{k=0}^p (-1)^k (c_k-b_k)  \geq 0$ }

By Euler-Poincar\'e, this is zero for the entire sum. It appears as if the 
Witten deformatin proof (see e.g. \cite{Cycon}) works in the discrete too. 

\section{Isospectral deformation}

\paragraph{}
If $d$ is the exterior derivative, the operator $D=d+d^*$ is the {\bf Dirac operator} of $G$. 
The Dirac operator $D$ admits an {\bf isospectral Lax deformations} $D' = [B,D]=BD-DB$, where
$B=d-d^* + \gamma i b$, if $D=d+d^*+b$. The parameter $\gamma$ is a tuning parameter.
For $\gamma=0$ the deformation stays real. For $\gamma \neq 0$, it is allowed to become
complex. The Dirac operator $D(t)$ defines for every $t$ an {\bf elliptic complex} $D:E \to F$
meaning that we have a splitting $D(t): E \to F$ such that McKean-Singer relation holds. 
\index{Dirac operator}
\index{Lax deformation}
\index{Isospectral deformation}
\index{Elliptic complex} 

\satz{The Lax system for the Dirac operator is integrable. }

\paragraph{}
The spectrum of $D(t)$ stays constant. Actually, $L=D(t)^2$ stays constant. 
\index{Lax system}

\paragraph{}
We have a deformation of the complex for which all classical geometry like the wave equation
stays the same because $L$ does not change. It is only the underlying $d$ which changes. The 
{\bf Connes formula} $\sup_{|Df|_{\infty}=1} |f(x)-f(y)|$ allows to re-interpret 
the isospectral deformation as a deformation of the metric. 
\index{Connes formula}
\index{Distance formula ! Connes}

\section{Trees and Forests}

\paragraph{}
Given a finite simple graph $G$, a {\bf rooted spanning tree} is a subgraph $H$ of $G$ 
which is a tree with the same vertex set together with a base point $x$. 
A {\bf rooted spanning forest} is a subgraph $H$ of $G$ which is a forest with the
same vertex set together with a base point $x$. Let $K$ be the {\bf Kirchhoff Laplacian} 
of the graph and ${\rm Det}(K)$ the {\bf pseudo determinant}, the product of the non-zero
eigenvalues of $K$. It is ${\rm exp}(-\zeta'(0))$ for the zeta function of $K$. 
\index{Spanning tree}
\index{Rooted tree}
\index{Rooted spanning tree}
\index{Pseudo determinant}
\index{Kirchhoff Laplacian}
\index{Laplacian Kirchhoff}

\paragraph{}
The {\bf tree number} of a graph $G$ is the number of rooted spanning tree in $G$. 
The {\bf forest number}  of a graph is the number of rooted spanning forests. 
The first part of the following theorem is the {\bf Kirchhoff matrix tree theorem}.
The second part of the theorem is the {\bf Chebotarev-Shamis forest theorem}. 
\index{Kirchhoff matrix tree theorem}
\index{Chebotarev-Shamis forst theorem}
\index{Forest theorem}
\index{Matrix tree theorem}

\satz{${\rm Det}(K)$ is the tree number. ${\rm det}(K+1)$ is the forest number.}

\index{Tree number}
\index{Forest number}

\paragraph{}
By Baker-Norine theory, the tree number is also the order of the {\bf Picard group} which
appears in the context of discrete Riemann-Roch. 
\index{Picard group}

\paragraph{}
If $F,G$ are arbitrary $n \times m$ matrices. 
Assume $p(x) = p_0 (-x)^m + p_1 (-x)^{m-1} + \cdots + p_k (-x)^{m-k} + \cdots + p_m$
is the {\bf characteristic polynomial} of the $m \times m$ matrix $F^T G$ with $p_0=1$.
The {\bf generalized Cauchy-Binet theorem} is 

\satz{ $p_k = \sum_{|P|=k} \det(F_P) \det(G_P)$ }

where the sum is over $k$-minors and where $p_k$ are the coefficients
of the characteristic polynomial of $F^T G$. It implies the
polynomial identity
$\det(1+z F^T G) = \sum_P z^{|P|} \det(F_P) \det(G_P)$
in which the sum is over all minors $A_P$ including the empty one $|P|=0$ for which
$\det(F_P) \det(G_P)=1$.

\section{Wave equation}

\paragraph{}
Because the Hodge Laplacian is a square $L=D^2 = (d+d^*)^2$, the wave equation $u_{tt} = L u$ has
an explicit d'Alembert solution. Let $D^{-1}$ be the pseudo inverse of $D$. It is
defined as $\sum_{k, \lambda_k \neq 0} u_k u_k^T/\lambda_k$, where $D u_k=\lambda_k u_k$
with an orthonormal eigenbasis $\{u_k\}$ of $D$. 
\index{d'Alembert solution}
\index{Hodge Laplacian}

\satz{
$u(t) = \cos(D t) u(0) + i \sin(D t) D^{-1} u'(0)$
}

\paragraph{}
With the {\bf complex wave} $\psi(t) = u(t) - i D u'(0)$, we can write 
the solution of the {\bf real wave} equation of $u$ as a solution of the 
{\bf Schr\"odinger equation}.
\index{Schroedinger equation}
\index{Wave equation}

\satz{
$\psi(t) = e^{i D t} \psi(0)$.
}

\paragraph{}
Just use the Euler identity $e^{i D t} = \cos(Dt) + i \sin(D t)$ and plug in 
$\psi(t) = u(t) - i D u'(0)$ to see that the relation holds. 

\section{Euler-Poincar\'e}

\paragraph{}
Let $\Lambda^p(G)$ be the functions from $G_p=\{ x \in G \; | \; {\rm dim}(x)=k \; \}$ to $R$
which are anti-symmetric. The {\bf exterior derivatives} 
$$   d_pf(x_0,x_1,\dots, x_p)=\sum_{j} (-1)^j f(x_0,\dots,\hat{x_j},\dots,x_p) $$
define linear map $d:\Lambda(G) \to \Lambda(G)$, where $\Lambda(G)$ is the Hilbert space
of dimension $n=|G|$. Since $d^2=0$, the {\bf cohomology groups} $H^p(G)={\rm ker}(d_p)/{\rm im}(d_{p-1})$ are
defined. Their dimensions are the Betti numbers $b_p(G)$. 
The matrix $H=(d+d^*)^2$ decomposes into blocks $H_k(G)$. We have the {\bf Hodge relations}: 
\index{Exterior derivative}
\index{Cohomology group}
\index{Hodge relation}

\satz{
${\rm dim}({\rm ker}(H_k)) = {\rm dim}(H^k)$. 
}

\paragraph{}
Define the {\bf Poincar\'e polynomial} $p_G(t) = \sum_{k=0} b_k(G) t^k$. The {\bf cohomological
Euler characteristic} is $p_G(-1) = b_0(G)-b_1(G)+b_2(G) - \cdots$. If the $f$-vector of $G$
is $(v_0,v_1,v_2, \dots)$, then the {\bf Euler polynomial} is $e_G(t) = \sum_{k=0} v_k(G) t^k$. 
By definition, we have $d_G(-1)=\chi(G)$. The {\bf Euler-Poincar\'e theorem} tells that the combinatorial
and cohomological Euler characteristic agree. 
\index{Euler-Poincare theorem}
\index{Euler polynomial}
\index{Poincar\'e polynomial}
\index{Polynomial ! Euler}
\index{Polynomial ! Poincar\'e}

\satz{ $\chi(G) = e_G(-1) = p_G(-1)$.  }

\section{Interaction cohomology}

\paragraph{}
Let $\Lambda^p_2(G)$ be the functions from 
$G_p = \{ (x,y) \; | \; x \cap y \neq \emptyset, {\rm dim}(x)+{\rm dim}(y)=p \}$
which are anti-symmetric. Like {\bf Stokes theorem} $df(x)=f(\delta x)$ for simplicial cohomology,
we define the {\bf exterior derivative} $df((x,y))=f(\delta{x},y)+ (-1)^{{\rm dim}}(x) f(x,\delta y)$
with the understanding that $f(\delta{x},y)=0$ if $\delta{x} \cap y=\emptyset$ or
$f(x,\delta{y})=0$ if $x \cap \delta{y}=\emptyset$. 
It defines a linear map $d:\Lambda_2(G) \to \Lambda_2(G)$, where $\Lambda_2(G)$ has
as dimension the number of intersecting simplices $(x,y)$ in $G$. Again, we can define 
the {\bf Dirac operator} $D=d+d^*$ and the {\bf Hodge operator} $H=D^2$ and decompose 
the later into blocks $H_k$. As before:

\index{Dirac operator ! Connection cohomology}
\index{Stokes theorem}
\index{Exterior derivative}

\satz{
${\rm dim}({\rm ker}(H_k)) = {\rm dim}(H^k)$.
}

\paragraph{}
The {\bf quadratic Poincar\'e polynomial} $p_G(t) = \sum_{k=0} b_k(G) t^k$
and {\bf quadratic Euler polynomial} $e_G(t) = \sum_{k=0} v_k(G) t^k$ are defined in the same way.
By definition, we have $d_G(-1)=\chi(G)$. The {\bf Euler-Poincar\'e theorem} tells that the combinatorial
and cohomological Wu characteristic agree.

\satz{ $\omega(G) = e_G(-1) = p_G(-1)$.  }

\index{Poincar\'e polynomial | Wu characteristic}
\index{Poincar\'e polynomial | quadratic}

\section{Stokes theorem}

\paragraph{}
Examples of {\bf orientation oblivious} measurements are valuations $F$ like
$F(A)=v_k(A)$ measuring the $k$ dimensional volume of a subcomplex $A$ of $G$
or $\chi(A)$ giving the Euler characteristic of a subcomplex. The {\bf length}
of a subcomplex $A$ for example is $v_1(A)$. In the continuum, such quantities
are accessible via {\bf integral geometry}, like {\bf Crofton type formulas}.
In the discrete one refers to it also as {]bf geometric probability theory}.
\index{Orientation oblivious} 
\index{Integral geometry}
\index{Geometric probability}
\index{Crofton}

\paragraph{}
If valuations are done after an orientation has been chosen on the elements of $G$,
we get a {\bf calculus} which features a {\bf fundamental theorem}. 
Given an arbitrary choice of orientation of the sets in $G$, the boundary $\delta A$
of a subcomplex is in general no more a subcomplex, it becomes a {\bf chain}. 
Given a form $F \in \Lambda$, we can still compute $F(\delta A)$. 
If $G$ is {\bf orientable} $d$-complex and $A$ is a $k$-subcomplex with {\bf boundary}
$\delta A$, then $\delta A$ is a complex. {\bf Stokes theorem} tells that
for any $k$-subcomplex $A$ with boundary $\delta A$, and any $k$-form $F$
\index{Calculus} 
\index{Fundamental theorem ! calculus}
\index{Orientable}
\index{Boundary} 

\satz{$dF(A) = F(\delta A)$.}

\paragraph{}
For $k=1$, we talk about the {\bf fundamental theorem of line integrals}, for $k=2$ we have {\bf Stokes theorem}
and $k=3$ goes under the name {\bf divergence theorem}. The derivative $d_0: \Lambda^0 \to \Lambda^1$ is the
{\bf gradient}, the derivative $d_1: \Lambda^1 \to \Lambda^2$ is the {\bf curl} and $d_2: \Lambda^2 \to \Lambda^3$ 
is the divergence (often just identified with the dual $d_0^*: \Lambda^1 \to \Lambda^0$, as $2$-forms and $1$-forms
in three dimensions are dual to each other). 
This Stokes theorem holds both for the familiar {\bf simplicial calculus} related to Euler characteristic $\chi(G)$
as well as the {\bf connection calculus} related to the Wu characteristics $\omega_k(G)$. 
\index{Connection calculus}
\index{Simplicial calculus} 
\index{Curl}
\index{Gradient}
\index{Divergence}

\section{Quadratic Lefschetz fixed point}

\paragraph{}
Given an automorphism $T$, 
define the {\bf quadratic Lefschetz number} $\chi_T(G)$, the super trace of the induced map on cohomology. 

\satz{
$\chi_T(G)=\sum_{x \sim y, (x,y)=(T(x),T(y))} i_T(x,y)$
}

\index{Lefschetz formula | Wu characteristic}

\paragraph{}
We can especially look at the case when $G$ is a ball. This is cohomologically non-trivial. 

\satz{
An endomormorphism of a ball $G$ has a fixed $(x,y)$, $x \cap y \neq \emptyset$. 
}

\index{Brouwer fixed point | quadratic} 

\section{Eulerian spheres}

\paragraph{}
Let $\G_d$ be the class of {\bf $d$-graphs}, $\Sphere_d$ the class of {\bf $d$-spheres},
$\B_d$ the class of {\bf $d$-balls}, and $\C_k$ the class of graphs with {\bf chromatic number }
$k$. Note that all Barycentric refinements of a complex are Eulerian. 
We call the class $\Sphere_d \cap \C_{d+1}$ the class of {\bf Eulerian spheres}
and $\B_d \cap C_{d+1}$ the class of Eulerian disks. The $0$-sphere $2$ is Eulerian.
Eulerian $1$-spheres are cyclic graphs with an even number of vertices. 
\index{Eulerian sphere}
\index{Chromatic number}
\index{Eulerian sphere}

\satz{
Every unit sphere of an Eulerian sphere is Eulerian.
}

\paragraph{}
The {\bf dual graph} $\hat{G}$ of a $d$-sphere $G$ is the graph in which 
the $d$-simplices are the vertices and where two simplices are connected, 
if one is contained in the other. A graph $(V,E)$ is {\bf bipartite} if $V=(A \cup B$ with disjoint
$A,B$ such $E=\{ (a,b) \; | \;  a \in A, b \in B \}$. Every Barycentric refinement of a complex
is a bipartite graph as we can take $A = \{ x \in G \; {\rm dim}(x) \; {\rm even} \}$ and 
$B =  \{ x \in G \; {\rm dim}(x) \; {\rm odd} \}$. 
\index{Dual graph} 
\index{Bipartite graph}

\satz{
For $G \in \Sphere_d$, then $\hat{G}$ is bipartite if and only if $G$ is Eulerian.
}

\section{Riemann-Hurwitz}

\paragraph{}
The {\bf automorphism group} ${\rm Aut}(G)$ of a simplicial complex is the group
of all {\bf automorphisms} of $G$. An {\bf endomorphism} $T$ is a simplicial map $G \to G$. 
If an endomorphism $T$ is restricted to the {\bf attractor} $\bigcap_k T^k(G)$ is
an automorphism. An automorphism $T$ of $G$ induces automorphisms on Barycentric 
refinements and so {\bf graph automorphisms}. The equivalence classes $G_1/A$ are
graphs.

\index{Automorphism group}
\index{Endomorphisms}
\index{Graph automorphism}
\index{Attractor of an automorphism}

\satz{ If $A \subset {\rm Aut}(G)$, then $G_1/A$ is a simplicial complex.}

\paragraph{}
We can see $G_1$ as a {\bf branched cover} $G_1/A$, {\bf ramified} over some points. 
If $G$ was a $d$-graph, then $G_1/A$ is a discrete {\bf orbifold}. If there
are no ramification points, then the cover $G \to G/A$ is a fibre bundle with
structure group $A$. 
\index{Orbifold}
\index{Branched cover}
\index{Ramified}

\paragraph{}
Given an automorphism $T$, define the {\bf ramification index} 
$e(x)=1-\sum_{T \neq 1, T(x)=x} \omega(x)$ of $X$. 
The following remark was obtained with  Tom Tucker. It is a discrete
{\bf Riemann-Hurwitz} result:
\index{Ramification index}
\index{Riemann-Hurwitz}

\satz{ $\chi(G) = |A| \chi(G/A) - \sum_{x \in G} (e(x)-1)$ }

\paragraph{}
For every subset $\G_k$ of indices of fixed dimension $k$, we have
by the {\bf Burnside lemma} $\sum_{T \in A} \sum_{x \in \G_k, T(x)=x} 1 = |A| |\G_k|$.
The super sum gives $\sum_{T \in A} \sum_{x, T(x)=x} \omega(x) = |A| \chi(H)$.
This gives 
$\sum_{T \neq 1} \sum_{x \in G} \omega(x) + \sum_{x \in G} \omega(x)  = |A| \chi(H)$. 
\index{Burnside lemma}

\paragraph{}
Let $\chi(G,T)$ denote the {\bf Lefschetz number} of $T$. From the Lefschetz fixed point
formula we get
\index{Lefschetz number}

\satz{
$\chi(G/A) = \frac{1}{|A|} \sum_{T \in A} L(G,T)$ 
}

\section{Riemann-Roch}

\paragraph{}
A {\bf divisor} $X$ is an integer-valued function on $G$. The {\bf simplex Laplacian}
$L$ is defined as $L(x,y)=\omega(x) \omega(y) H_0(x,y)$, where $H_0$ is 
the Kirchhoff Laplacian of the {\bf simplex graph} in which $G$ is the vertex
set and two $x,y$ are connected if one is contained in the other and the 
dimensions differ by $1$. The simplex graph is one-dimensional as it has no triangles. 
A divisor $X$ is called {\bf principal} if $X=Lf$ for some integer valued function $f$. 
We think of a divisor as a geometric object and define the Euler characteristic 
$\chi(G)=\sum_x \omega(x) X(x)$. A divisor is {\bf essential} if $\omega(x) X(x) \geq 0$ 
for all $x$. The {\bf linear system} $|X|$ of $X$ is the set of
$f$ for which $X+(f)$ is essential. Its {\bf dimension} $l(X)$ is the 
maximal $k \geq 0$ such that for every $m<k$ and every $Y$ of $\chi(Y)=m$, the divisor
$X-Y$ is essential. Define the {\bf canonical divisor} $K(x)=0$. The simplest Riemann-Roch theorem is 
\index{Divisor}
\index{Simplex Laplacian}
\index{Simplex graph}
\index{Laplacian ! Simplex}
\index{Graph ! Simplex}
\index{Linear system}
\index{Principal divisor}
\index{Essential divisor}
\index{Divisor ! essential}
\index{Divisor ! principal}

\satz{ $l(X) - l(K-X) = \chi(X)$. }

\paragraph{} 
This is {\bf Baker-Norine theory}, slightly adapted to change the perspective: classically a divisors appear
one a one dimensional connected curve (Riemann surface or 1-dimensional graph) $G$ and ${\rm deg}(X)+\chi(G) = \chi(G)$. 
Centering at the geometric underlying object gives the canonical divisor $K=-2$ which is in the case when $G$ is one-dimensional is 
linearly equivalent to the negated curvature function $K(v)=-2+{\rm deg}(v)$ on the {\bf vertices} of $G$. 
Riemann-Roch tells that the signed distance to the surface $\chi(G)=0$ is $\chi(G)$. 
\index{Riemann Roch}
\index{Riemann Surface}
\index{Baker-Norine}

\paragraph{}
Reflecting at $0$ rather than at usual canonical divisor representing the curve $G$
allows to have a Riemann-Roch for arbitrary dimensions. Generalizing Baker-Norine naively to higher dimensional 
simplicial complexes does not work, as the curvature $\kappa$ of $\chi(G)$ has only in the one-dimensional
case the property that $K=-2\kappa$ is a divisor. Classically $l(X), L(K-X)$ have cohomological interpretations.
Also here, Riemann-Roch appears like a fancy Euler-Poincar\'e formula, but it is deeper than the later,
as surface ${\rm ker}(\chi)$ is {\bf bumpy}: it contains both {\bf generic divisors} as well as {\bf special divisors}.
\index{Generic divisor}
\index{Special divisor}

\paragraph{}
The image of $L$ is a linear subspace of the set ${\rm ker}(G) = \chi(G)=0$. The quotient ${\rm ker}(\chi)/{\rm im}(L)$ 
is the {\bf Picard group} or {\bf divisor class group}. The equivalence classes of divisors can be represented
by rooted spanning trees in the simplex graph. This defines a group structure on {\bf rooted spanning trees}.
That there is a bijective identification between divisor classes and spanning trees is the subject of: 

\index{Divisor class group}
\index{Picard group}
\index{Rooted spanning trees}

\satz{ The Picard group is isomorphic to the tree group.}

\chapter{References}

\paragraph{}
For the history of topology\cite{Dieudonne1989,HistoryTopology} 
and graph theory \cite{BiggsLloydWilson,HistoryTopology,handbookgraph} 
and discrete geometry \cite{BobenkoSuris}. See \cite{Hatcher,Spanier,Rotman} for
notations in algebraic topology, \cite{HararyGraphTheory,Biggs,BM} for graph theory.

\paragraph{}
Abstract simplicial complexes appeared in 1907 by Dehn and Heegaard
\cite{BurdeZieschang,MunkholmMunkholm}. In \cite{alexandroff} they appeared under the name
{\bf unrestricted skeleton complex}. In \cite{Whitehead}, J.H.C. Whitehead calls them
{\bf symbolic complexes}. 
\index{Symbolic complex}
\index{Complex ! Symbolic}
\index{Whitehead}
\index{Dehn}
\index{Heegaard}
\index{Skeleton complex}
\index{Complex ! Skeleton}

\paragraph{}
Some of the results generalize to $\Delta$ sets or simplicial sets. Some connection
calculus however does not. Some connection calculus does not go over yet. The
unimodularity theorem does not hold for simplicial sets, at least for the approaches we 
tried so far. 

\paragraph{}
Homotopy theory as developed by \cite{Whitehead} uses elementary
expansions and contractions. Homotoptic complexes are said to have
the same ``nucleus". \cite{Whitehead} uses ``collapsible" for ``homotopic to a point".
See also \cite{WhiteheadI}. The notions appearing for simplices described
by graph theory, see \cite{I94a,I94,CYY}. 
\index{Homotopy}
\index{Whitehead}
\index{Elementary expansions}

\paragraph{}
Dimension theory has a long history \cite{Crilly}. The inductive definition of graphs
appeared first in \cite{elemente11}. We studied the average in \cite{randomgraph}.
\index{Inductive dimension}

\paragraph{}
{\bf Random graphs} were first studied in \cite{erdoesrenyi59}. The average Euler characteristic
appears in  \cite{randomgraph}. 
\index{Random graphs}

\paragraph{}  
The idea of seeing geometric quantities as expectations is central in {\bf integral geometry}.
The first time, that curvature was seen as an expectation of indices is Banchoff 
\cite{Banchoff1967,Banchoff1970}. Random methods in geometry is part of integral geometry as pioneered 
by Crofton and Blaschke \cite{Blaschke,Nicoalescu}.
We have used in in \cite{colorcurvature, indexexpectation} and
\cite{indexformula}. Having curvature given as an expectation allows to deform it. Given a
unitary flow $U_t$ on functions for example produces a deformation of the curvature.
\index{Integral geometry}
\index{Banchoff}
\index{Crofton}
\index{Blaschke}

\paragraph{}
Discrete curvature traces back to a combinatorial curvature considered by Heesch \cite{Heesch} in the
context of graph coloring and extended in \cite{Gromov87}.
The formula $K(p)=1-V_1(p)/6$ and for graphs on the sphere appears also in
\cite{Presnov1990,Presnov1991}, where it is also pointed out that $\sum_p K(p)=2$
is Gauss-Bonnet formula. Discrete curvature was used in \cite{Higuchi} and unpublished work of Ishida from 1990.
Higushi use $K(p) = 1-\sum_{y \in S(p)} (1/2-1/d(y))$, where $d(y)$ are the cardinalities of 
the neighboring face degrees in the sphere $S(p)$. For two dimensional graphs, where all faces are triangles,
this simplifies to $d_j=3$ so that $K=1-|S|/6$, where $|S|$ is the cardinality of the
sphere $S(p)$. In \cite{elemente11} second order curvatures were used.
The {\bf Levitt curvature} in arbitrary dimension appears in \cite{Levitt1992}. 
We rediscovered it in \cite{cherngaussbonnet} 
after tackling dimension by dimension separately, not aware of Levitt. We got into the topic 
while working on \cite{elemente11}. Chern's proof is \cite{Chern44} followed 
\cite{AllendoerferWeil,Fenchel}.
See \cite{Rosenberg, Cycon} for modern proofs.  Historical remarks are in \cite{Chern1990}.
\index{Heesch}
\index{Ishida} 
\index{Higushi}
\index{Chern}
\index{Levitt curvature}

\paragraph{}
The Erd\"os R\'enyi probability space were introduced in \cite{erdoesrenyi59}. 
The formulas for the average dimension and Euler characteristic has been found in 
\cite{randomgraph}. The recursive dimension was first used in \cite{elemente11}. 
We looked at more functionals in \cite{KnillFunctional}. 
\index{Erdoes-Renyi}

\paragraph{}
The {\bf discrete Hadwiger Theorem} appears in \cite{KlainRota}. The continuous version is \cite{Hadwiger}.
For integral geometry and geometric probability, see \cite{Santalo}.
The theory of valuations on distributive lattices has been pioneered by Klee \cite{Klee63}
and Rota \cite{Rota71} who proved that there is a unique
valuation such that $X(x)=1$ for any join-irreducible element. See also \cite{forman2000}. 
\index{Klee}
\index{Rota}
\index{Diescrete Hadwiger}

\paragraph{}
Wu characteristic appeared in \cite{Wu1953} and was discussed in \cite{Gruenbaum1970}. 
We worked on it in \cite{valuation} and 
announced cohomology in \cite{CaseStudy2016} and \cite{MathTableMarch}. 
For the {\bf connection cohomology} belonging to Wu characteristic, see \cite{CohomologyWuCharacteristic}. 
\index{Connection Cohomology}
\index{Cohomology Connection}

\paragraph{}
For discrete Poincar\'e-Hopf see \cite{poincarehopf} and an attempt to popularize it in 
\cite{knillcalculus} or Mathematica demonstrations \cite{KnillWolframDemo1,KnillWolframDemo2}.
It got pushed a bit more in \cite{indexformula}. 
For the classical Poincar\'e-Hopf, see \cite{Spivak5}. 
For the classical case, Poincar\'e covered the 
$2$-dimensional case in chapter VIII of \cite{poincare85}
It got extended by Hopf in arbitrary dimensions \cite{hopf26}. It is pivotal in the 
proof of Gauss-Bonnet theorems for smooth Riemannian manifolds 
(i.e. \cite{Guillemin,Spivak1999,Hirsch,Henle1994,docarmo94,BergerGostiaux}).

\paragraph{}
Discrete McKean-Singer was covered in \cite{knillmckeansinger}. 
The best proof in the continuum is \cite{Cycon}. The classical result is 
\cite{McKeanSinger}. 
In \cite{DiscreteAtiyahSingerBott}, the suggestion appeared to define elliptic
discrete complexes using McKean-Singer. 
\index{McKean-Singer}

\paragraph{}
The Zykov sum (join) was introduced in \cite{Zykov} to graph theory. 
The strong ring was covered in \cite{ArithmeticGraphs,StrongRing}. 
\index{Zykov sum}

\paragraph{}
The Brouwer-Lefschetz theorem is \cite{brouwergraph}. It generalizes
the $1$-dimensional case \cite{NowakowskiRival}.
The classical result is \cite{Lefschetz26a}. See also \cite{Hopf28}.
\index{Brouwer-Lefschetz}

\paragraph{}
The classical {\bf Kuenneth formula} is \cite{Kuenneth}. The graph version \cite{KnillKuenneth},
uses the Barycentric refinement $(A \times B)_1$ of the Cartesian product $A \times B$. 
\index{Kuenneth formula}

\paragraph{}
About the history of discrete notions of manifolds, see \cite{Scholz}. The Evako definition
of a sphere as a cell complex for which every unit sphere is a $n-1$ sphere and such that
removing one point makes it contractible was predated by approaches 
of Vietoris or van Kampen. The later would have accepted homology spheres as unit spheres.
\index{Evako Sphere}
\index{Vietoris}
\index{Van Kampen}
\index{Homology sphere}

\paragraph{}
The classical Sard theorem is \cite{Sard42}. The discrete version was remarked in \cite{KnillSard}. 
\index{Sard Theorem}

\paragraph{}
For the spectral universality, see \cite{KnillBarycentric} and \cite{KnillBarycentric}.
It uses a result of Lidskii-Last \cite{SimonTrace} which 
assures if $||\mu-\lambda||_1 \leq \sum_{i,j=1}^{n} |A-B|_{ij}$
for any two symmetric $n \times n$ matrices $A,B$ with
eigenvalues $\alpha_1 \leq \alpha_2 \leq \dots \leq \alpha_n$. 
\index{Lidskii-Last}

\paragraph{}
The discrete exterior derivative goes back to Betti and Poincar\'e and was already anticipated
by Kirchhoff. As pointed out in \cite{brouwergraph}, the discrete Hodge point is \cite{Eckmann1}. 
It appeared also in \cite{DanijelaJost}. 
The {\bf discrete Dirac operator} was stressed in \cite{KnillILAS}. 
\index{Betti}
\index{Kirchhoff}
\index{Discrete Dirac operator}
\index{Dirac operator}

\paragraph{}
The unimodularity theorem $|{\rm det}(L)|=1$ was discovered in February 2016,
announced in \cite{MathTableOctober} and proven in \cite{Unimodularity}. 
An other proof was given in \cite{MukherjeeBera2018}. 
\index{Unimodularity}

\paragraph{}
We have looked at the arithmetic of unit spheres in \cite{Spheregeometry}, especially
in the context of the diagonal Green function entries. The other Green function entries
are covered in \cite{ListeningCohomology}.

\paragraph{}
The result $\chi(G)=p(G)-n(G)$ was proven in \cite{HearingEulerCharacteristic,ListeningCohomology}. 
The functional equation for the spectral zeta function of the connection Laplacian
was proven in \cite{DyadicRiemann}. Earlier work in the Hodge Zeta case is \cite{KnillZeta}.
The zeta function is called Dyadic because the Barycentric limit is in an ergodic setup a 
{\bf von Neumann-Kakutani system} \cite{Kni95}, which has the Pr\"ufer group as the spectrum. The system
is a group translation on the dyadic group of integers and also known as the {\bf adding machine}. 
\index{Dyadic zeta function} 
\index{Von Neumann-Kakutani} 
\index{Adding machine}
\index{Dyadic group}
\index{Pr\"ufer group}

\paragraph{}
The {\bf Hydrogen relation} $H=L-L^{-1}$ for one-dimensional 
complexes was studied in \cite{DehnSommerville,Helmholtz} and \cite{Hydrogen}.
\index{Hydrogen relation}

\paragraph{}
An earlier talk \cite{classicalstructures} summarizes things also.
\cite{DiracKnill} is an earlier snapshot about the linear algebra part. 
\cite{ArchimedesFunctions,knillcalculus} summarize the calculus. 

\paragraph{}
The matrix tree theorem is \cite{Kirc}. It is based on the {\bf Cauchy-Binet theorem}
\cite{Cauchy,Binet}. A generalization \cite{CauchyBinetKnill} gives the coefficients
of the characteristic polynomial. The {\bf Chebotarev-Shamis theorem} is 
\cite{ChebotarevShamis1,ChebotarevShamis2}. 
See also \cite{knillforest}, where we initially were not aware of the work of Chebotarev and Shamis.
\index{Chebotarev-Shamis}
\index{Matrix tree}
\index{Caucy-Binet}

\paragraph{}
The {\bf Lax deformation} of exterior derivatives was introduced in 
\cite{IsospectralDirac,IsospectralDirac2} and was motivated by {\bf Witten deformation}
\cite{Witten1982,Cycon}. Lax systems were introduced first to
\cite{Lax1968}. Commutation relations of that form have appeared earlier when 
describing {\bf free tops} $L'=[B,L]$, where $B=I^{-1} L$ is the angular velocity and $L$
the angular velocity in $so(n)$, which are geodesics in $SO(n)$ \cite{Arnold1980}.
\index{Free top}
\index{Geodesics}
\index{Lax deformation}
\index{Witten deformation}
\index{Isospectral deformation}

\paragraph{}
The Connes formula \cite{Connes} is elementary but crucial in the process of 
generalizing Riemannian geometry to {\bf non-commutative geometry}.
\index{Connes formula}
\index{Non-commutative geometry}

\paragraph{}
After finding a multiplication completing the Zykov addition to a ring 
in \cite{ArithmeticGraphs}, we realized it is the dual to the Sabidussi ring. 
In \cite{StrongRing}, we looked at the ring generated by the Cartesian product.
It is a subring and consists of discrete CW complexes. Unlike for simplicial sets,
the classical theorems like Gauss-Bonnet and energy theorem go over. 
\index{Strong ring}
\index{Arithmetic of graphs}
\index{Discrete CW complex}
\index{Zykov addition}

\paragraph{}
Riemann-Roch for graphs is \cite{BakerNorine2007}. See also \cite{BakerNorine}.
We worked on Riemann-Hurwitz in \cite{TuckerKnill}. The usual approach for
Riemann-Hurwitz in graph theory is to see them as discrete analogues of
algebraic curves or Riemann surfaces see \cite{MednykhMednykh}.
\index{Riemann-Hurwitz}
\index{Riemann-Roch}

\paragraph{}
\cite{Weyl1925} first looked for a combinatorial definition of spheres. Forman
\cite{forman95} defined spheres through the Reeb as objects admitting 2 critical points.
See also \cite{forman98}. More on discrete Morse theory in \cite{Forman1999, Forman2003}. 
\index{Reeb}
\index{Forman}
\index{Discrete Morse theory}

\paragraph{}
We used data fitting to get first heuristically the Stirling formula then proved it.
It is however considered "well known" \cite{BrentiWelker}. 
It appears also in \cite{Stanley86,LuzonMoron,Hetyei}. 

\paragraph{}
The history of polytopes is a ``delicate task" \cite{Devadoss}.
The Euler polyhedron formula (Euler's gem) was discussed in \cite{Richeson}. 
The early proofs of Schl\"afli and Staudt had still gaps according to 
\cite{BurdeZieschang}. The difficulty is also explained in \cite{lakatos,Gruenbaum2003}.

\paragraph{}
The story of polyhedra is told in \cite{Richeson,coxeter}. Historically, it was
developed in \cite{Schlafli}, \cite{Schoute}, \cite{AliciaBooleStott}.
Coxeter \cite{coxeter} defines a polytop as a convex body with polygonal faces. 
\cite{gruenbaum} also works with convex polytopes in $R^n$ where the dimension is the 
dimension of the affine span.

\paragraph{}
The perils of a general definition of a polytop were known since Poincar\'e 
(see \cite{Aczel,Richeson,cromwell,lakatos}).
Polytop definitions are given in \cite{Schlafli,coxeter,gruenbaum,symmetries}.
Topologists started with new definitions \cite{alexandroff,Fomenko,ConwayCapstone,Spanier},
and define first a simplicial complex and then polyhedra as topological spaces which admit a 
{\bf triangularization} by a simplicial complex.
\index{Convex polytop}
\index{Polytop}

\paragraph{}
Dehn-Sommerville relations have traditionally been formulated for
convex polytopes and then been generalized to situations where unit spheres
can be realized as convex polytopes. See
\cite{Klee1964,NovikSwartz,MuraiNovik,LuzonMoron,BrentiWelker,Hetyei,Klain2002}
or \cite{BergerLadder}. 
\index{Dehn-Sommerville}

\paragraph{}
We started to think about graph coloring during the project \cite{KnillNitishinskaya}. 
The reports \cite{knillgraphcoloring} and \cite{knillgraphcoloring2} explored this a bit more. 
It is related to Fisk theory \cite{Fisk1977b,Fisk1977a}.
\index{Fisk theory}

\paragraph{}
Some special graphs appearing when counting was considered in \cite{CountingAndCohomology}.
When writing this, we were not aware that the cell complex introduced already in \cite{Bjoerner2011}
which goes much further than what we did. Other classes of complexes called {\bf orbital networks} 
 \cite{KnillOrbital1,KnillOrbital2,KnillOrbital3} were studied first with Montasser Ghachem. 
\index{Orbital networks}
\index{Bjoerner complex}
\index{Ghachem}

\paragraph{}
For the Alexander duality, see \cite{BjoernerTancer}. Originally established by
Alexander in 1922, it was formulated by Kalai and Stanley in combinatorial
topology. We formulated it with cohomology rather than homology and cohomology. As such
it is an identity where we have numbers on both sides. 
\index{Alexander Duality}
\index{Kalai}
\index{Stanley}

\chapter{Questions}

\section{Inverse spectral questions}

\paragraph{}
We have seen that the spectrum of $L$ does not determine the Betti numbers in general
but that for a Barycentric refinement of $G$, the Betti numbers $b_0,b_1$ can 
be read of from the spectrum as the number of eigenvalues $1$ and $-1$. 

\question{Does the spectrum of $L$ determine $b_k$ for $k \geq 2$.}

\question{Does the spectrum of $L$ determine the Wu characteristic $\omega(G)$? }

\section{Barycentric limit}

We have seen that the limiting spectral measure can be computed in the case $d=1$.
It is a smooth measure. In higher dimensions, we see spectral gaps. These gaps have first been
seen in the {\bf BeKeNePaPeTe paper} \cite{hexacarpet}.
\index{BeKeNePaPeTe paper}

\question{Prove spectral gaps in limiting spectral measure for $d \geq 2$. }

\section{Coloring}

\paragraph{}
The {\bf four color theorem} is equivalent to the statement that all 2-spheres are 4-colorable. 
\index{Four color theorem}

\question{Are all $d$-spheres $(d+2)$-colorable?}

\question{Are all $2$-graphs $5$ colorable?}

\section{Connection Cohomology}

\paragraph{}
While we know that connection cohomology is not a homotopy invariant, we have not yet
proven that it is a topological invariant. We have introduced a notion of homeomorphism
in \cite{KnillTopology}. One can also use the notion whether geometric realizations are
homeomorphic to ask:

\question{ Is connection cohomology a topological invariant?}

\paragraph{}
We would like to find more examples of triangulations of non-homeomorphic d-manifolds with different
connection cohomology which can not be distinguished by other means:

\question{ Can one distinguish {\bf homology spheres} with Wu cohomology?}

\index{Homology sphere} 

\paragraph{}
Something we have only started to look at"

\question{ Is there a duality for connection cohomology? }

\paragraph{}
As connection cohomology is not a homotopy invariant, the naive generalization
does not work. 

\section{Random complexes}

\paragraph{}
The probability spaces $E(n,p)$ of graphs define natural random spaces of 
simplicial complexes as we can take the Whitney complex of a graph. While
we have a formula for the expectation of Euler characteristic, this is not
yet available for Wu characteristic numbers $\omega_k$. 

\question{ What is the expected value of $\omega_k$ on $E(n,p)$? }

\paragraph{}
We would also like to know the expectations of the Betti numbers:

\question{ What is the expectation of $b_k(G)$ on $E(n,p)$? }

\section{Zeta function}

\paragraph{}
While various equivalent expressions exist for the 
{\bf connection zeta function} in the Barycentric limit
of a one-dimensional complex, we don't yet have found a reference
about where the roots of $\zeta$ are:

\question{The limiting zeta function $\zeta$ has roots on the imaginary axes.} 

\printindex
\chapter{Bibliography}

\bibliographystyle{plain}

\end{document}